\documentclass[11pt]{amsart}
\usepackage{amsmath, amsfonts, upref}
\newtheorem{prop}{Proposition}[section]

\newtheorem{lemma}[prop]{Lemma}
\newtheorem{cor}[prop]{Corollary}
\newtheorem{theorem}[prop]{Theorem}

\begin{document}
\newcommand{\ee}{\mathbb{E}}
\newcommand{\pp}{\mathbb{P}}
\newcommand{\ii}{\mathbb{I}}
\newcommand{\zz}{\mathbb{Z}}
\newcommand{\rr}{\mathbb{R}}
\newcommand{\ex}{\mathrm{Exp}(1)}
\newcommand{\ra}{\rightarrow}
\newcommand{\fp}{f^\prime}
\newcommand{\fpp}{f^{\prime\prime}}
\newcommand{\ww}{W^\prime}
\newcommand{\frc}[2]{{\textstyle{\frac{#1}{#2}}}}
\newtheorem{thm}{Theorem}
\newtheorem{lmm}{Lemma}

\begin{center}
{\bf Stein's Method and Minimum Parsimony Distance after Shuffles}
\end{center}

\begin{center}
{\bf Running head: Stein's Method and Minimum Parsimony}
\end{center}

\begin{center}
Version of 10/25/04 
\end{center}

\begin{center}
By Jason Fulman (fulman@math.pitt.edu)
\end{center}
\begin{center}
University of Pittsburgh Math Department, 414 Thackeray Hall
\end{center}
\begin{center}
Pittsburgh, PA 15260
\end{center}

{\bf Abstract}: Motivated by Bourque and Pevzner's simulation study of
the parsimony method for studying genome rearrangement, Berestycki and
Durrett used techniques from random graph theory to prove that the
minimum parsimony distance after iterating the random transposition
shuffle undergoes a transition from Poisson to normal behavior. This
paper establishes an analogous result for minimum parsimony distance
after iterates of riffle shuffles or iterates of riffle shuffles and
cuts. The analysis is elegant and uses different tools: Stein's method
and generating functions. A useful technique which emerges is that of
making a problem more tractable by adding extra symmetry, then using
Stein's method to exploit the symmetry in the modified problem, and
from this deducing information about the original problem.

\begin{center}
\end{center}

\begin{center}
\end{center}

\begin{center}
\end{center}

\section{Introduction}

	In the study of genome rearrangement, one often views
genomes as signed permutations, where each integer corresponds to a
unique gene/marker and the sign corresponds to its orientation. For
unichromosomal genomes, the most frequent rearrangements are
reversals. A reversal $\rho(i,j)$ applied to a permutation $\pi =
\pi_1 \cdots \pi_{i-1} \pi_i \cdots \pi_j \pi_{j+1} \cdots \pi_n$
reverses the segment $\pi_i \cdots \pi_j$ to obtain a new permutation
$\pi_1 \cdots \pi_{i-1} \ -\pi_j \ -\pi_{j-1} \cdots -\pi_i \pi_{j+1}
\cdots \pi_n$. For instance the reversal $\rho(2,4)$ would send $4 \ 1
\ -5 \ -2 \ 3$ to $4 \ 2 \ 5 \ -1 \ 3$. In this context, the minimum
parsimony distance of a signed permutation $\pi$ is defined as the
minimum number of reversals needed to bring the identity permutation
to $\pi$. Hannenhalli and Pevzner \cite{HP} found an exact
combinatorial formula for the minimum parsimony distance of a signed
permutation; for this and other results, see the book \cite{P}.

	A fundamental problem, emphasized in \cite{BP} and the survey
	\cite{Du}, is to understand the distribution of the minimum
	parsimony distance after a given number of reversals has
	occurred. More generally for any shuffling technique on
	permutations (which from now on we assume are unsigned so that
	we are dealing with the symmetric group $S_n$), one can define
	the minimum parsimony distance of a permutation $\pi$ as the
	number of shuffles needed to bring the identity to $\pi$. And
	it is natural to study the distribution of the minimum
	parsimony distance after a given number of shuffles has
	occurred.

	An exciting recent work along this lines is the paper
	\cite{BeDu}, which studied the random transposition walk on the
	symmetric group on n symbols. Then the minimum parsimony
	distance of $\pi$ is simply n $-$ number of cycles of
	$\pi$. Letting $D_t$ be the minimum parsimony distance after t
	iterations of the random transposition walk, they showed that
	$D_{cn/2} \sim u(c)n$, where $u$ is an explicit function
	satisfying $u(c)=c/2$ for $c \leq 1$ and $u(c)<c/2$ for
	$c>1$. They also described the fluctuation of $D_{cn/2}$ about
	its mean in each of three regimes (subcritical where the
	fluctuations are Poisson, critical, and supercritical where
	the fluctuations are normal). They exploit a connection
	between the transposition walk and random graphs (about which
	an enormous amount is known).

	In the current paper we examine minimum parsimony distance for
a more vigorous shuffling method, the Gilbert-Shannon-Reeds riffle
shuffle. While we do not know if this is of biological interest, the
mathematical ubiquity of riffle shuffles (see the survey \cite{Di} for
an overview of connections to dynamical systems, Lie theory and much
else) as well as possible applications in casinos more than justifies
the question. Riffle shuffling proceeds as follows. Given a deck of
$n$ cards, one cuts it into 2 packets with probability of respective
pile sizes $j,n-j$ given by $\frac{{n \choose j}}{2^n}$. Then cards
are dropped from the packets with probability proportional to the
packet size at a given time; thus if the current packet sizes are
$A_1,A_2$, the next card is dropped from packet $i$ with probability
$A_i/(A_1+A_2)$. Bayer and Diaconis \cite{BD} prove the fundamental
result that after $r$ riffle shuffles, the probability of obtaining
the permutation $\pi^{-1}$ is $\frac{{n+2^r - d(\pi) - 1 \choose
n}}{2^{rn}}$. Here $d(\pi)$ denotes the number of descents of $\pi$,
that is $|\{i: 1 \leq i \leq n-1, \pi(i)>\pi(i+1) \}|$. For instance
the permutation $3 \ 1 \ 4 \ 2 \ 5$ has two descents. From the
Bayer-Diaconis formula it is clear that the minimum parsimony distance
of a permutation $\pi^{-1}$ is simply $\lceil log_2(d(\pi)+1)
\rceil$. Thus the study of minimum parsimony distance for riffle
shuffles is the study of the distribution of $d(\pi)$ under the
measure $\frac{{n+2^r - d(\pi) - 1 \choose n}}{2^{rn}}$.

	More generally, for k and n integers, we let $R_{k,n}$ denote
	the measure on $S_n$ which chooses $\pi$ with probability
	$\frac{{n+k-d(\pi)-1 \choose n}}{k^n}$ and study the number of
	descents. First let us review what is known. As $k \rightarrow
	\infty$, the distribution $R_{k,n}$ tends to the uniform
	distribution on $S_n$. It is well known (\cite{CKSS}),
	\cite{T}) that for $n \geq 2$ the number of descents has mean
	$\frac{n-1}{2}$ and variance $\frac{n+1}{12}$ and that
	$\frac{d(\pi)-(n-1)/2}{\sqrt{(n+1)/12}}$ is asymptotically
	normal. Aldous \cite{A} proved that $\frac{3}{2} log_2(n)$
	riffle shuffles are necessary and suffice to be close to the uniform distribution on $S_n$. Bayer and Diaconis \cite{BD} give more refined
	asymptotics, proving that for $k=2^c n^{3/2}$ with $c$ a real
	number,
\[ \frac{1}{2}
 \sum_{\pi \in S_n} |R_{k,n}(\pi)-\frac{1}{n!}| = 1 - 2
 \Phi(\frac{-2^{-c}}{4 \sqrt{3}}) + O(n^{-1/4}) \] where \[ \Phi(x) =
 \frac{1}{\sqrt{2 \pi}} \int_{- \infty}^x e^{-t^2/2} dt.\] Motivated by this result, Mann \cite{Ma}
 proved that if $k=a n^{3/2}$, with $a$ fixed, then the number of
 descents has mean $\frac{n-1}{2} - \frac{n^{1/2}}{12 a} + O(1)$ and variance
 $\frac{n}{12} + O(n^{1/2})$, and is asymptotically normal as $n
 \rightarrow \infty$. He deduces this from Tanny's local limit theorem
 for $d(\pi)$ under the uniform distribution \cite{T} and from the formula for
 $R_{k,n}$.

	We prove two new results concerning the distribution of
	$d(\pi)$ under the measure $R_{k,n}$. First, we complement the above results on normal approximation by using Stein's
	method to upper bound the total
	variation distance between the distribution of $k-1-d(\pi)$
	and a Poisson variable with mean $\frac{k}{n+1}$; our bound shows
	the approximation to be good when $\frac{k}{n}$ is small. Second, we use generating functions to give very
	precise asymptotic estimates for the mean and variance of
	$d(\pi)$ when $\frac{k}{n}>\frac{1}{2 \pi}$. For instance if $k=\alpha
	n$ with $\alpha>\frac{1}{2 \pi}$, we show that
\[ \left| \ee_{R_{k,n-1}}(d+1) - n \left( \alpha -
  \frac{1}{e^{1/\alpha}-1} \right) - \left( \frac{e^{1/\alpha} ( -2
\alpha e^{1/\alpha} + 2\alpha + e^{1/\alpha}+1)}{2 \alpha^2
(e^{1/\alpha}-1)^3} \right) \right|\] is at most
$\frac{C_{\alpha}}{n}$ and that \[ \left| Var_{R_{k,n-1}}(d) - n
\left( \frac{e^{1/\alpha} (\alpha^2 e^{2/\alpha} + \alpha^2 - 2
\alpha^2 e^{1/\alpha} - e^{1/\alpha})}{\alpha^2 (e^{1/\alpha}-1)^4}
\right) \right| \leq A_{\alpha} \] where $C_{\alpha}, A_{\alpha}$ are
constants depending on $\alpha$ (and are independent of $\alpha$ for
$\alpha \geq 1$). Mann \cite{Ma} had exact expressions for the mean
and variance (which we derive another way) but only obtained
asymptotics when $k=a n^{3/2}$ with $a$ fixed. The main technical
point and effort of the paper \cite{SGO} on information loss in card
shuffling was to obtain asymptotics for the mean and variance of
$d(\pi)$ under the measure $R_{k,n}$. By a clever application of the
method of indicator variables, they obtain bounds, but ours are much
better.

	Next let us describe the technique we use to study the
	distribution of $d(\pi)$ under $R_{k,n}$, as we believe this
	to be as interesting as the result itself. To apply Stein's
	method to study a statistic $W$, one often uses an
	exchangeable pair $(W,W')$ of random variables  (this
	means that the distribution of $(W,W')$ is the same as that of
	$(W',W)$) such that the conditional expectation $\ee(W'|W)$ is
	approximately $(1-\lambda)W$. Typically to construct such a
	pair one would use a Markov chain on $S_n$ which is reversible
	with respect to the measure $R_{k,n}$, choose $\pi$ from
	$R_{k,n}$, let $\pi'$ be obtained from $\pi$ from one step
	in the chain, and finally set $(W,W')=(W(\pi),W(\pi'))$. For
	the problem at hand this does not seem easy. Thus we modify
	the problem; instead of considering the measure $R_{k,n}$, we
	consider a measure $C_{k,n}$ which chooses a permutation $\pi$
	with probability $\frac{{n+k-c(\pi)-1 \choose n-1}}{n
	k^{n-1}}$. Here $c(\pi)$ is the number of cyclic descents of
	$\pi$, defined as $d(\pi)$ if $\pi(n)<\pi(1)$ and as
	$d(\pi)+1$ if $\pi(n)>\pi(1)$. This probability measure was
	introduced in \cite{F1} and $C_{2^r,n}(\pi)$ gives the chance
	of obtaining $\pi^{-1}$ after first cutting the deck at a
	uniformly chosen random position and then performing r
	iterations of a riffle shuffle. The advantage of working with $C_{k,n}$ is that it has a natural
	symmetry which leaves it invariant, since performing two
	consecutive cuts at random positions is the same as performing
	a single cut. As will be explained in Section \ref{exch}, this
	symmetry leads to an exchangeable pair $(d,d')$ with the very
	convenient property that $\ee_C(d'|\pi)$ is approximately
	$(1-\frac{1}{n})d$. We obtain a Poisson approximation theorem
	for $d$ under the measure $C_{k,n}$. Although the measures
	$R_{k,n}$ and $C_{k,n}$ are not close when $\frac{k}{n}$ is
	small (a main result of \cite{F3} is that that the total
	variation distance between them is roughly $\frac{n}{4k}$ for
	$k \geq n$), we show that the distribution of $k-d$ under
	$C_{k,n}$ is close to the distribution of $k-c$ under
	$C_{k,n}$ which in turn is equal to the distribution of
	$k-1-d$ under $R_{k,n-1}$. This implies a Poisson
	approximation theorem for the original problem of interest.

	Incidentally, it is proved in \cite{F1} that r iterations of
``cut and then riffle shuffle" yields exactly the same distribution as
performing a single cut and then iterating r riffle shuffles. Thus the
chance of $\pi^{-1}$ after r iterations of ``cut and then riffle
shuffle" is $C_{2^r,n}(\pi)$, which implies that the minimum parsimony
distance of the ``cut and then riffle shuffle" process is $\lceil
log_2(c(\pi)) \rceil$. Thus the study of minimum parsimony distance
for the ``cut and then riffle shuffle" procedure is equivalent to the
study of $c$ under the distribution $C_{k,n}$. But as mentioned in the
last paragraph, we will prove that this is the same as the
distribution of $d+1$ under $R_{k,n-1}$. Hence the theory of minimum
parsimony distance for ``cut and then riffle shuffle" is equivalent to
the theory for riffle shuffles, and we shall say nothing more about
it.

	The reader may wonder why we don't apply our exchangeable pair
	for normal approximation. Most theorems for Stein's method for
	normal approximation assume that the pair $(W,W')$ satisfies
	the property $\ee(W'|W)=(1-\lambda) W$ for some $\lambda$. In
	our case this only approximately holds, that is
	$\ee(W'|W)=(1-\lambda) W + G(W)$ where $G(W)$ is small. There
	are normal approximation results in the literature (\cite{RR},
	\cite{C}) for dealing with this situation, but they require
	that $\frac{\ee(|G(W)|)}{\lambda}$ goes to 0. Using interesting
	properties of Eulerian numbers, we show that even for the
	uniform distribution (the $k \rightarrow \infty$ limit of
	$C_{k,n}$) the quantity $\frac{\ee(|G(W)|)}{\lambda}$ is
	bounded away from 0. Finding a version of Stein's method which
	allows normal approximation for our exchangeable pair (even
	for the uniform distribution) is an important open
	problem. Incidentally, for the special case of the uniform
	distribution, it is possible to prove a central limit theorem
	for $d$ by Stein's method \cite{F2}, using a different
	exchangeable pair.

	Having described the main motivations and ideas of the paper,
	we describe its organization. Section \ref{exch} defines an
	exchangeable pair to be used in the study of $d(\pi)$ under
	the measure $C_{k,n}$, and develops a number of properties of
	it. It also gives closed formulas (but not asymptotics) for
	the mean and variance of $d(\pi)$, by relating them to the
	mean and variance of $c(\pi)$, and computing the latter using
	generating functions. Section \ref{poisson} uses the
	exchangeable pair of Section \ref{exch} to prove a Poisson
	approximation theorem for $k-d(\pi)$ under the measure
	$C_{k,n}$ (the bounds are valid for all integer values of $k$
	and $n$ but informative only when $\frac{k}{n}$ is small). It
	then shows how to deduce from this a Poisson approximation
	theorem for $k-1-d(\pi)$ under the measure $R_{k,n}$. Section
	\ref{other} gives asymptotics for the mean and variance for
	$c(\pi)$ under $C_{k,n}$ for $\frac{k}{n} > \frac{1}{2 \pi}$
	(and so also for $d(\pi)$ under $C_{k,n}$ and $R_{k,n}$). It
	then explores further properties of the exchangeable pair
	which are related to normal approximation. Finally, it gives a
	quick algorithm for sampling from $R_{k,n}$, which should be
	useful in empirically investigating the nature of the
	transition from Poisson to normal behavior.

\section{The Exchangeable pair, mean, and variance} \label{exch}

	This section constructs an exchangeable pair $(d,d')$ for the
	measure $C_{k,n}$ and develops some of its
	properties. Throughout we let $\ee_C$ denote expectation with
	respect to $C_{k,n}$. We relate $\ee_C(d)$ and $\ee_C(d^2)$ to
	$\ee_C(c)$ and $\ee_C(c^2)$, and then use generating functions
	to find expressions (whose asymptotics will be studied later)
	for $\ee_C(c)$ and $\ee_C(c^2)$.

	To begin let us construct an exchangeable pair $(d,d')$. We
	represent permutations $\pi$ in two line form. Thus the
	permutation represented by
\[ \begin{array}{c c c c c c c c c}
                i \ & : & 1 & 2 & 3 & 4 & 5 & 6 & 7 \\ \pi(i) & : & 6
                & 4 & 1 & 5 & 3 & 2 & 7 \end{array} \] sends 1 to 6, 2
                to 4, and so on. One constructs a permutation $\pi'$
                by choosing uniformly at random one of the $n$ cyclic
                shifts of the symbols in the bottow row of the two
                line form of $\pi$. For instance with probability 1/7
                one obtains the permutation $\pi'$ which is
                represented by \[ \begin{array}{c c c c c c c c c} i \
                & : & 1 & 2 & 3 & 4 & 5 & 6 & 7 \\ \pi'(i) & : & 5 & 3
                & 2 & 7 & 6 & 4 & 1 \end{array} .\] An essential point
                is that if $\pi$ is chosen from the measure $C_{k,n}$,
                then so is $\pi'$; note that this would not be so for
                the measure $R_{k,n}$. Thus if one chooses $\pi$ from
                $C_{k,n}$, defines $\pi'$ as above, and sets
                $(d,d')=(d(\pi),d(\pi'))$, it follows that $(d,d')$ is
                exchangeable with respect to the measure
                $C_{k,n}$. Observe also that $d'-d \in \{0, \pm 1\}$.

	Recall that $\pi$ is said to have a cyclic descent at position
	$j$ if either $1 \leq j \leq n-1$ and $\pi(j)>\pi(j+1)$ or
	$j=n$ and $\pi(n)>\pi(1)$. It is helpful to define random
	variables $\chi_j(\pi)$ ($1 \leq j \leq n$) where
	$\chi_j(\pi)=1$ if $\pi$ has a cyclic descent at position $j$
	and $\chi_j(\pi)=0$ if $\pi$ does not have a cyclic descent at
	position $j$. We let $\ii$ denote the indicator function of an
	event. We also use the standard notion that if $Y$ is a random
	variable, $\ee(Y|A)$ is the conditional expectation of $Y$
	given $A$.

\begin{lemma} \label{steinmean}
 \[ \ee_C(d'-d|\pi) = -\frac{d}{n} + \frac{n-1}{n}
\ii_{\chi_n(\pi)=1}. \] \end{lemma} \begin{proof} Note that $d'=d+1$
occurs only if $\pi$ has a cyclic descent at $n$ and that then it
occurs with probability $\frac{n-1-d}{n}$. Note also that $d'=d-1$
occurs only if $\pi$ does not have a cyclic descent at $n$, and that
it then occurs with probability $\frac{d}{n}$. To summarize,
\begin{eqnarray*}
\ee_C(d'-d|\pi) & = & -\frac{d}{n} \ii_{\chi_n(\pi)=0} + \frac{n-1-d}{n} \ii_{\chi_n(\pi)=1}\\
& = & -\frac{d}{n} + \frac{n-1}{n} \ii_{\chi_n(\pi)=1}. \end{eqnarray*} \end{proof}

	As a corollary, we obtain $\ee_C(d)$ in terms of $\ee_C(c)$.

\begin{cor} \label{cormean} \[ \ee_C(d) = \frac{n-1}{n} \ee_C(c).\]
\end{cor} \begin{proof} Since $(d,d')$ is an exchangeable pair,
$\ee_C(d'-d)=0$. It follows that $\ee_C(\ee_C(d'-d|\pi))=0$. So from Lemma
\ref{steinmean} \[ \ee_C(d) = (n-1) \ee_C(\ii_{\chi_n(\pi)=1}).\] Since
the variables $\chi_1(\pi),\cdots,\chi_n(\pi)$ have the same
distribution under $C_{k,n}$, and $c=\chi_1(\pi)+\cdots+\chi_n(\pi)$,
the result follows. \end{proof}

	Lemma \ref{use1} will be helpful at several points in this paper.

 \begin{lemma} \label{use1} \[ \ee_C(d \ii_{\chi_n(\pi)=1}) =
\ee_C \left( \frac{c(c-1)}{n} \right) .\] \end{lemma}
 
\begin{proof} Observe that 
\begin{eqnarray*}
\ee_C(d \ii_{\chi_n(\pi)=1}) & = & \ee_C((c-1) \ii_{\chi_n(\pi)=1})\\
& = & \frac{1}{n} \sum_{i=1}^n \ee_C((c-1) \ii_{\chi_i(\pi)=1})\\
& = & \frac{1}{n} \ee_C(c(c-1)). \end{eqnarray*} \end{proof}

	As a consequence, we obtain $\ee_C(d^2)$ in terms of $\ee_C(c)$ and
	$\ee_C(c^2)$.

\begin{cor} \label{2ways} \[ \ee_C(d^2) = (1-\frac{2}{n}) \ee_C(c^2) + \frac{1}{n} \ee_C(c).\] \end{cor}

\begin{proof} \begin{eqnarray*} \ee_C(c^2) & = & \ee_C(d+\ii_{\chi_n(\pi)=1})^2\\
& = & \ee_C(d^2) + 2 \ee_C(d \ii_{\chi_n(\pi)=1}) + \ee_C(\ii_{\chi_n(\pi)=1})\\
& = & \ee_C(d^2) + \frac{2}{n} \ee_C (c^2-c) + \frac{1}{n} \ee_C(c), \end{eqnarray*} where the final equality is Lemma \ref{use1}. This is equivalent to the statement of the corollary. \end{proof}

	Next we use generating functions to compute $\ee_C(c)$ and
	$\ee_C(c^2)$. For this some lemmas are useful.

\begin{lemma} \label{ful} (\cite{F1}) For $n>1$, the number of elements in $S_n$ with $i$ cyclic
 descents is equal to n multiplied by the number of elements in
 $S_{n-1}$ with $i-1$ descents.
\end{lemma} 

\begin{lemma} \label{genfunc} \[ \frac{\sum_{\pi \in S_n} t^{c(\pi)}}{(1-t)^n} = n \sum_{m \geq 0} m^{n-1} t^m .\] \end{lemma}

\begin{proof} Given Lemma \ref{ful}, the result now follows from the well
 known generating function for descents (e.g. \cite{FS}) \[
 \frac{\sum_{\pi \in S_n} t^{d(\pi)+1}}{(1-t)^{n+1}} = \sum_{m \geq 0}
m^n t^m.\]
\end{proof}

	Proposition \ref{exactmean} gives a closed formula for
	$\ee_C(c)$.

\begin{prop} \label{exactmean} For $n>1,$ \[ \ee_C(c) = k - \frac{n}{k^{n-1}} \sum_{j=1}^{k-1} j^{n-1} .\]
\end{prop}

\begin{proof} Multiplying
 the equation of Lemma \ref{genfunc} by $(1-t)^n$ and then
 differentiating with respect to $t$, one obtains the equation is \[
 \sum_{\pi \in S_n} c(\pi) t^{c(\pi)-1} = n(1-t)^n \sum_{m \geq 0}
 m^nt^{m-1} - n^2 (1-t)^{n-1} \sum_{m \geq 0} m^{n-1}t^{m}.\]
 Multiplying both sides by $\frac{t}{nk^{n-1} (1-t)^n}$ gives the
 equation \[ \sum_{\pi \in S_n} \frac{c(\pi) t^{c(\pi)}}{nk^{n-1}
 (1-t)^n} = \frac{1}{k^{n-1}} \sum_{m \geq 0} m^nt^{m} -
 \frac{n}{k^{n-1} (1-t)} \sum_{m \geq 0} m^{n-1}t^{m+1}.\] The
 coefficient of $t^k$ on the left hand side is precisely the expected
 value of $c$ under the measure $C_{k,n}$. The proposition now
 follows by computing the coefficient of $t^k$ on the right hand side.
\end{proof}

	By a similar argument, one obtains an exact expression for
	$\ee_C(c^2)$.

\begin{prop} \label{exact2} For $n>1$,
\[ \ee_C(c^2)  = k^2 - \frac{n(n+1)}{k^{n-1}} \sum_{j=1}^{k-1} j^n + \frac{n(nk-n-k)}{k^{n-1}} \sum_{j=1}^{k-1} j^{n-1}.\]
\end{prop}

\begin{proof} From the proof of Proposition \ref{exactmean}, we know that \[ 
 \sum_{\pi \in S_n} c(\pi) t^{c(\pi)} = n(1-t)^n \sum_{m \geq 0}
 m^nt^{m} - n^2 (1-t)^{n-1} \sum_{m \geq 0} m^{n-1}t^{m+1}. \]
 Differentiate with respect to $t$, multiply both sides by
 $\frac{t}{nk^{n-1} (1-t)^n}$, and take the coefficient of $t^k$. On the left hand side one gets $\ee_C(c^2)$. On the right hand side one obtains the
 coefficient of $t^k$ in \begin{eqnarray*} & & \frac{1}{k^{n-1}} \sum_{m \geq 0}
 m^{n+1}t^m - \frac{2n}{k^{n-1} (1-t)} \sum_{m \geq 0} m^n t^{m+1}\\ & & -
 \frac{n}{k^{n-1} (1-t)} \sum_{m \geq 0} m^{n-1}t^{m+1} +
 \frac{n(n-1)}{k^{n-1}(1-t)^2} \sum_{m \geq 0} m^{n-1}t^{m+2}.\end{eqnarray*} After elementary simplifications the
 result follows. \end{proof}

\section{Poisson regime} \label{poisson}

	A main result of this section is a Stein's method proof that for $k$ much smaller
	than $n$, the random variable $X(\pi):=k-d(\pi)$ under the
	measure $C_{k,n}$ is approximately Poisson with mean
	$\lambda:=\frac{k}{n}$. Then we show how
	this can be used to deduce Poisson limits for $k-c(\pi)$ under
	the measure $C_{k,n}$ and for $k-1-d(\pi)$ under the measure $R_{k,n}$.

	To begin we recall Stein's method for Poisson approximation. A
	book length treatment of Stein's method for Poisson
	approximation is \cite{BHJ}, but that book emphasizes the
	coupling approach. We prefer to work from first principles
	along the lines of Stein's original formulation as presented
	in \cite{Stn}.
 
	Throughout we use the exchangeable pair $(X,X')$, where $\pi$
and $\pi'$ are as in Section \ref{exch}, $X(\pi)=k-d(\pi)$,
$X'=X(\pi')$, and the underlying probability measure is $C_{k,n}$. Let
$\pp_{\lambda}$ denote probability under the Poisson distribution of
mean $\lambda$, and as usual let $\pp_C$ denote probability with
respect to the measure $C_{k,n}$. Let $A$ be subset of $\zz^+$, the
set of non-negative integers. Stein's method is based on the following
``Stein's equation'' \[ \pp_C(X \in A) - \pp_{\lambda} \{A \} = \ee_C
(i T_{\lambda} - T \alpha) g_{\lambda,A}.\]

	 Let us specify the terms on the right hand side of the
equation, in the special case of interest to us.

\begin{enumerate}
\item The function $g=g_{\lambda,A}: \zz^+ \mapsto \rr$ is constructed
\ to solve the equation
\[ \lambda g(j+1) - jg(j) = \ii_{j \in A} - \pp_{\lambda} \{A\} , \ j \geq 0 \] where $g(0)$ is taken to be 0. We also need the following lemma which bounds certain quantities related to $g$.

\begin{lemma} \label{supg} (\cite{BHJ}, Lemma 1.1.1)
\begin{enumerate}
\item Let $||g||$ denote $sup_j |g_{\lambda,A}(j)|$. Then $||g|| \leq 1$ for all $A$. 
\item Let $\Delta(g)$ denote $sup_j |g_{\lambda,A}(j+1)-g_{\lambda,A}(j)|$. Then $\Delta(g) \leq 1$ for all $A$.
\end{enumerate}
\end{lemma}

\item The map $T_{\lambda}$ sends real valued functions on
$\zz^+$ to  real valued functions on
$\zz^+$ and is defined by $T_{\lambda}(f)[j] = \lambda
f(j+1)-j f(j)$.

\item The map $i$ sends real valued functions on $\zz^+$ to real
valued functions on $S_n$, the symmetric group. It is defined by
$(if)[\pi] = f (X(\pi))$.

\item The map $T$ is a map from the set of real valued antisymmetric
functions on $S_n \times S_n$ to the set of real valued functions on
$S_n$. It is defined by $Tf[\pi] = \ee_C(f(\pi,\pi')|\pi)$. (Since the
pair $(\pi,\pi')$ is exchangeable, $\ee_C(Tf)=0$, which is crucial for
the proof of the Stein equation).

\item Finally (and this is where one has to make a careful choice),
the map $\alpha$ is a map from real valued functions on $Z^+$ to
antisymmetric real valued functions on $S_n \times S_n$. In the proof of Theorem
\ref{poisapprox} we will specify which $\alpha$ we use.

\end{enumerate}

	In order to approximate $X$ by a Poisson($\lambda$) random
	variable, it will be useful to approximate the mean of the
	random variable $\frac{c(\pi)}{n}$ by $\lambda$. This is
	accomplished in the next lemma, the second part of which is
	not needed in the sequel.

\begin{lemma} \label{approx}  Let $\lambda=\frac{k}{n}$, where $k,n$ are positive integers.
\begin{enumerate}
\item \[ \left| \frac{\ee_C(c)}{n} - \lambda \right| \leq k(1-\frac{1}{k})^n.\]
\item  \[ \left| \frac{\ee_C(c)}{n} - \lambda \right| \leq \frac{k}{n}.\]
\end{enumerate}
\end{lemma}

\begin{proof} For the
first assertion, note that by Proposition \ref{exactmean}, \[
\left|\frac{\ee_C(c)}{n} - \lambda \right| = \frac{1}{k^{n-1}}
\sum_{j=1}^{k-1} j^{n-1} \leq \frac{(k-1)^n}{k^{n-1}}.\] The second
assertion follows since the formula for $C_{k,n}$ forces $c \leq k$
with probability 1. \end{proof}

\begin{theorem}
 \label{poisapprox} Let $\lambda=\frac{k}{n}$ where $k,n$ are positive
 integers. Then for any $A \subseteq \zz^+$, \[ |\pp_C(k-d(\pi) \in A) -
 \pp_{\lambda}(A)| \leq (\frac{k}{n})^2 + k(n+1)
 (1-\frac{1}{k})^n .\]
\end{theorem}

\begin{proof} As above, let $X(\pi)=k-d(\pi)$ and $X'(\pi) = k-d(\pi')$. Since $A,\lambda$ are fixed, throughout the
 proof the function $g_{\lambda,A}$ is denoted by $g$. We specify the
 map $\alpha$ to be used in the ``Stein equation'' \[ \pp_C(X \in A) -
 \pp_{\lambda} \{A \} = E_C (i T_{\lambda} - T \alpha) g.\] Given a real
 valued function $f$ on $\zz^+$, we define $\alpha f$ by \[ \alpha f
 [\pi_1,\pi_2] = f(X(\pi_2)) \ii_{X(\pi_2)=X(\pi_1)+1} - f(X(\pi_1))
 \ii_{X(\pi_1)=X(\pi_2)+1} .\] Note that as required this is an antisymmetric
 function on $S_n \times S_n$.

	Then one computes that $T \alpha g$ is the function on $S_n$
	defined by \[ T \alpha g(\pi) = \ee_C(\alpha g(\pi,\pi')|\pi) = \ee_C
	\left(g(X') \ii_{X'=X+1} - g(X) \ii_{X=X'+1}|\pi \right) .\]
	Thus by the reasoning of Lemma \ref{steinmean}, \[ T \alpha g (\pi) =
	  g(X+1) \ii_{\chi_n(\pi)=0}
	\frac{c(\pi)}{n} - g(X) \ii_{\chi_n(\pi)=1}
	(1-\frac{c(\pi)}{n}).\] Since $iT_{\lambda}g (\pi) = \lambda g(X+1) - X g(X)$ and
$\ii_{\chi_n(\pi)=0} = 1-\ii_{\chi_n(\pi)=1},$ one concludes that
\begin{eqnarray*} & & (i T_{\lambda}g - T \alpha g) [\pi] \\ & = &
\left[ (\lambda-\frac{c(\pi)}{n}) g(X+1) \right] + \left[ (\ii_{\chi_n(\pi)=1} - X) g(X)   \right]\\ & & +  \left[
\frac{c(\pi)}{n} \ii_{\chi_n(\pi)=1} (g(X+1)-g(X)) \right]\\
& = & \left[ (\lambda-\frac{c(\pi)}{n}) g(X+1) \right] + \left[ (c(\pi) - k) g(X)   \right]\\ & & +  \left[
\frac{c(\pi)}{n} \ii_{\chi_n(\pi)=1} (g(X+1)-g(X)) \right]. \end{eqnarray*}

	Thus to complete the proof, for each of the three terms in
	square brackets, we bound the expectation under the measure
	$C_{k,n}$. Lemma
	\ref{supg} and part 1 of Lemma \ref{approx} give that
	\begin{eqnarray*} \left| \ee_C \left( g(X+1)(\lambda-\frac{c(\pi)}{n}) \right) \right| &
	\leq & ||g|| \ee_C(|\lambda-\frac{c(\pi)}{n}|) \\ &
	\leq &  \ee_C(|\lambda-\frac{c(\pi)}{n}|)    \\ & \leq & k
	(1-\frac{1}{k})^n. \end{eqnarray*} For the second term in
	square brackets, one argues as for the first term in square
	brackets to get an upper bound of $nk (1-\frac{1}{k})^n$.

 To bound the expectation of the
	third term in square brackets, note that the nonnegativity of
	$c(\pi) \ii_{\chi_n(\pi)=1}$ and Lemma \ref{supg} imply that
	\begin{eqnarray*} & &\left| \ee_C \left[ \frac{c(\pi)}{n}
	\ii_{\chi_n(\pi)=1} \left(g(X+1)-g(X) \right) \right] \right|\\
	& \leq & \Delta(g) \ee_C (\frac{c(\pi)}{n}
	\ii_{\chi_n(\pi)=1})\\ & \leq & \ee_C (\frac{c(\pi)}{n}
	\ii_{\chi_n(\pi)=1})\\ & = & \ee_C(\frac{c(\pi)}{n^2} \sum_{i=1}^n \ii_{\chi_i(\pi)=1})\\ & = & \ee_C( \frac{c(\pi)^2}{n^2}). \end{eqnarray*} By the explicit formula for $C_{k,n}$, it follows that $c(\pi) \leq k$ with probability 1. Hence this is at most $(\frac{k}{n})^2$. \end{proof}

	To conclude this section, we show how Theorem \ref{poisapprox}
	can be used to deduce Poisson approximations for the two
	statistics we really care about: $k-c(\pi)$ under the measure
	$C_{k,n}$ and $k-1-d(\pi)$ under the measure $R_{k,n}$.

\begin{prop} \label{1sttransfer} For all $A \subseteq \zz^+$, \[ |\pp_C(k-d(\pi)
 \in A) - \pp_C(k-c(\pi) \in A) | \leq \frac{2k}{n}.\] \end{prop}

\begin{proof} Observe that for any $l \geq 0$,
\begin{eqnarray*}
& & \pp_C(d=l)\\ & = & \pp_C(d'=l)\\ & = & \pp_C(d'=l,c=l) +
\pp_C(d'=l,c=l+1)\\ & = & \pp_C(c=l) \pp_C(d'=l|c=l) + \pp_C(c=l+1)
\pp_C(d'=l|c=l+1)\\ & = & \pp_C(c=l) \frac{n-l}{n} + \pp_C(c=l+1)
\frac{l+1}{n}.
\end{eqnarray*} Thus
 \[ \pp_C(d=l) - \pp_C(c=l) = -\frac{l}{n} \pp_C(c=l) + \frac{l+1}{n}
 \pp_C(c=l+1) \] which implies that \[ |\pp_C(d=l) - \pp_C(c=l)| \leq 
 \frac{l}{n} \pp_C(c=l) + \frac{l+1}{n} \pp_C(c=l+1). \] Summing over $l \geq 0$ gives that \begin{eqnarray*} \sum_{l \geq 0} |\pp_C(d=l) - \pp_C(c=l)| & \leq & \sum_{l \geq 0} \left( \frac{l}{n} \pp_C(c=l) + \frac{l+1}{n} \pp_C(c=l+1) \right)\\
& = & \frac{2}{n} \sum_{l \geq 0} l \pp_C(c=l)\\
& = & \frac{2}{n} \ee_C(c)\\
& \leq & \frac{2k}{n}, \end{eqnarray*} where the final inequality is Proposition \ref{exactmean} (or also since $c \leq k$ with probability 1). The result follows. \end{proof}

\begin{cor} \label{cor1} Let $\lambda=\frac{k}{n}$ where $k,n$ are positive
 integers. Then for any $A \subseteq \zz^+$, \[ |\pp_C(k-c(\pi) \in A)
 - \pp_{\lambda}(A)| \leq (\frac{k}{n})^2 + \frac{2k}{n} + k(n+1)
 (1-\frac{1}{k})^n .\]
\end{cor}

\begin{proof} This is immediate  from Theorem
	\ref{poisapprox} and Proposition
	\ref{1sttransfer}. \end{proof}

	Proposition \ref{2ndtransfer} shows that the distribution of
	$d(\pi)+1$ under the measure $R_{k,n}$ is exactly the same
	as the distribution of $c(\pi)$ under the measure
	$C_{k,n+1}$.

\begin{prop} \label{2ndtransfer} For any $r \geq 0$, \[ \pp_{R_{k,n}}(d(\pi)=r) = \pp_{C_{k,n+1}}(c(\pi)=r+1).\] \end{prop}

\begin{proof} Note
 from the formula for $R_{k,n}$ that for any $r$, the probability of
 r descents under the measure $R_{k,n}$ is $\frac{{n+k-r-1 \choose
 n}}{k^n}$ multiplied by the number of permutations in $S_n$ with $r$
 descents. Similarly, from the formula for $C_{k,n+1}$ one sees that
 for any $r$, the probability of $r+1$ cyclic descents under the measure
 $C_{k,n+1}$ is $\frac{{n+k-r-1 \choose n}}{(n+1)k^n}$ multiplied by
 the number of permutations in $S_{n+1}$ with $r+1$ cyclic descents. The
 result follows from Lemma \ref{ful}. \end{proof} 

\begin{cor}  Let $\lambda=\frac{k}{n+1}$ where $k,n$ are positive
 integers. Then for any $A \subseteq \zz^+$, \begin{eqnarray*} && |\pp_{R_{k,n}}
(k-1-d(\pi) \in A) - \pp_{\lambda}(A)|\\& \leq & (\frac{k}{n+1})^2 +
\frac{2k}{n+1} + k(n+2) (1-\frac{1}{k})^{n+1} .\end{eqnarray*}
\end{cor}

\begin{proof} This is immediate Corollary \ref{cor1} and Proposition \ref{2ndtransfer}. \end{proof}

\section{Other regimes} \label{other}

	This section is organized into three subsections. Subsection
	\ref{asy} gives good asymptotic results for the mean and
	variance of $c(\pi)$ under $C_{k,n}$ (and so also for $d(\pi)$
	under $D_{k,n}$). Subsection \ref{further} develops further
	properties of the exchangeable pair $d,d'$ under
	$C_{k,n}$ which are relevant to normal approximation. Subsection \ref{sample} gives a speedy algorithm
	for sampling from the measures $C_{k,n}$ and $R_{k,n}$.

\subsection{Asymptotics of mean and variance} \label{asy}

	This subsection derives sharp estimates for the mean and
	variance of $c(\pi)$ under the measure $C_{k,n}$ when
	$\frac{k}{n} \geq \frac{1}{2 \pi}$. Since by Proposition
	\ref{2ndtransfer} the distribution of $d(\pi)$ under $R_{k,n}$
	is the same as the distribution of $c(\pi)-1$ under
	$C_{k,n+1}$, one immediately obtains results (which we stated in the introduction) for the mean and variance of $d(\pi)$ under
	$R_{k,n}$. We also remark that Corollaries \ref{cormean} and
	\ref{2ways} imply results for $d(\pi)$ under $C_{k,n}$.

	Throughout we will use information about the Bernoulli numbers
	$B_n$. They are defined by the generating function $f(z) =
	\sum_{n \geq 0} \frac{B_n z^n}{n!} = \frac{z}{e^z-1}$ so that
	$B_0=1,B_1=-\frac{1}{2},B_2=\frac{1}{6},B_3=0,B_4=-\frac{1}{30}$
	and $B_i=0$ if $i \geq 3$ is odd. The zero at 0 in the
	denominator of $f(z)$ cancels with the zero of $z$, so $f(z)$
	is analytic for $|z| < 2 \pi$ but has first order poles at $z=
	\pm 2 \pi i, \pm 4 \pi i, \cdots$. We also use the notation
	that $(n)_t$ denotes $n(n-1)\cdots(n-t+1)$ and that $(n)_0=1$.

	To see the connection with Bernoulli numbers, Lemma
	\ref{exactmean2} shows how to write $\ee_C(c)$ in terms of
	them.

\begin{lemma} \label{exactmean2} \[ \ee_C(c) = - k \sum_{t=1}^{n-1} \frac{B_t (n)_t}{t! k^t} .\]
\end{lemma} 

\begin{proof} This
 follows from Proposition \ref{exactmean} and by the expansion of
 partial power sums in \cite{GR}: \[ \sum_{r=0}^{a-1} r^n =
 \frac{a^n}{n+1} \left( a + \sum_{t=0}^{n-1} \frac{B_{t+1}
 (n+1)_{t+1}}{(t+1)! a^t} \right).\] \end{proof}

	Lemmas \ref{estimate0} and \ref{estimate1} give two elementary
	estimates.

\begin{lemma} \label{estimate0} For $0 \leq t \leq n$, \[ \left| \left(1-\frac{(n)_t}{n^t} \right) - \frac{{t \choose 2}}{n} \right| \leq \frac{{t \choose 2}^2}{2n^2} .\] \end{lemma}

\begin{proof} For $t=0,1$ the result is clear, so suppose that $t \geq 2$. We show that \[ \frac{{t \choose 2}}{n} -
	\frac{{t \choose 2}^2}{2n^2} \leq 1-\frac{(n)_t}{n^t} \leq
	\frac{{t \choose 2}}{n}.\] To see this write
	$1-\frac{(n)_t}{n^t} = 1-(1-\frac{1}{n})(1- \frac{2}{n}) \cdots
	(1-\frac{t-1}{n})$. One proves by induction that if
	$0<x_1,\cdots,x_n<1$, then $\prod_{j} (1-x_j) \geq 1 - \sum
	x_j$. Thus \[ (1-\frac{1}{n}) \cdots (1-\frac{t-1}{n}) \geq
	1-\sum_{j=1}^{t-1} \frac{j}{n} = 1-\frac{{t \choose 2}}{n} \]
	which proves the upper bound. For the lower bound note that
	\begin{eqnarray*} (1-\frac{1}{n}) \cdots (1-\frac{t-1}{n}) & =
	& e^{log(1-1/n) \cdots (1-(t-1)/n)}\\ & \leq &
	e^{-(\frac{1}{n}+\cdots+\frac{(t-1)}{n})}\\ & = & e^{-\frac{{t \choose
	2}}{n}}\\ & \leq & 1 - \frac{{t \choose 2}}{n} + \frac{{t
	\choose 2}^2}{2n^2}. \end{eqnarray*} The last inequality is on page 103 of \cite{HLP}. \end{proof}

\begin{lemma} \label{estimate1}  Suppose that $\alpha > \frac{1}{2 \pi}$. Then for $n \geq 1$, \[ \sum_{t=n}^{\infty}
\frac{|B_t| t^l}{\alpha^t t!}  \leq \frac{C_{l,\alpha} n^l}{(2 \pi
\alpha)^n} \] where $C_{l,\alpha}$ is a constant depending on $l$ and
$\alpha$ (and if $\alpha \geq 1$ the constant depends only on $l$). \end{lemma}

\begin{proof} Recall that $B_t$ vanishes for $t \geq 3$ odd
and that there is a bound $|B_{2t}| \leq 8 \sqrt{\pi t} (\frac{t}{\pi
e})^{2t}$ for $t \geq 1$ \cite{Le}.  Combining this with  Stirling's
bound $t! \geq  \sqrt{2 \pi t} (\frac{t}{e})^t e^{1/(12t+1)}$
\cite{Fe} one concludes that \[ \sum_{t=n}^{\infty} \frac{|B_t|
t^l}{\alpha^t t!} \leq C \sum_{t=n}^{\infty} t^l (2 \pi \alpha)^{-t} = \frac{C}{(2 \pi \alpha)^n} \sum_{t=0}^{\infty} \frac{(t+n)^l}{(2 \pi \alpha)^{t}} \leq \frac{C n^l}{(2 \pi \alpha)^n} \sum_{t=0}^{\infty} \frac{(1+t)^l}{(2 \pi \alpha)^t}
\] where $C$ is a universal constant. The ratio test shows
 that $ \sum_{t=0}^{\infty} \frac{(1+t)^l}{(2 \pi \alpha)^t} $ converges for $2
\pi \alpha >1$ (and moreover is at most a constant depending on $l$ if
$\alpha \geq 1$). \end{proof}

	We also require a result about Bernoulli numbers.

\begin{lemma} \label{bernoulli} Suppose that $\alpha>\frac{1}{2 \pi}$. Then

\begin{enumerate}
\item $\sum_{t=0}^{\infty} \frac{B_t}{t! \alpha^t} =
\frac{1}{\alpha(e^{1/\alpha}-1)}$.

\item $\sum_{t=0}^{\infty} \frac{B_t {t \choose 2}}{t! \alpha^t}=
\frac{e^{1/\alpha} ( -2 \alpha e^{1/\alpha} + 2\alpha +
e^{1/\alpha}+1)}{2 \alpha^3 (e^{1/\alpha}-1)^3}$.

\item $\sum_{t=0}^{\infty} \frac{B_{t+1}}{t! \alpha^t} = \frac{\alpha e^{1/\alpha} - e^{1/\alpha} - \alpha}{\alpha (e^{1/\alpha}-1)^2}$.

\item $\sum_{t=0}^{\infty} \frac{B_{t+1} {t \choose 2}}{t! \alpha^t} = \frac{e^{1/\alpha}(3 \alpha e^{2/\alpha} - e^{2/\alpha} -4 e^{1/\alpha}- 3\alpha -1)}{2 \alpha^3 (e^{1/\alpha}-1)^4}$.  \end{enumerate}
\end{lemma}

\begin{proof} We use the generating function $f(z)= \sum_{t
 \geq 0} \frac{B_t z^t}{t!} = \frac{z}{e^z-1}$ for the Bernoulli
 numbers, which as mentioned earlier is analytic for $|z| < 2
 \pi$. For the first assertion simply set $z=1/\alpha$. For the second
 assertion, one computes $\frac{ z^2}{2} \frac{d}{dz}
 \frac{d}{dz} f(z)$ and evaluates it at $z=1/\alpha$. For the third
 equation one differentiates $f(z)$ with respect to $z$ and then sets
 $z=\frac{1}{\alpha}$. For the fourth equation one differentiates
 $f(z)$ three times with respect to $z$, then multiplies by
 $\frac{z^2}{2}$ and sets $z=\frac{1}{\alpha}$. \end{proof}

	Next we give an estimate for $\ee_C(c)$.

\begin{prop} \label{critmean} Suppose
 that $k= \alpha n$ with $\alpha > \frac{1}{2 \pi}$. Then
\[
  \left| \ee_{C}(c(\pi)) - n \left( \alpha -
  \frac{1}{e^{1/\alpha}-1} \right) - \left( \frac{e^{1/\alpha} ( -2 \alpha e^{1/\alpha} + 2\alpha + e^{1/\alpha}+1)}{2 \alpha^2 (e^{1/\alpha}-1)^3} \right) \right| < \frac{C_{\alpha}}{n} \] where $C_{\alpha}$ is a constant depending on $\alpha$ (and which is independent of $\alpha$ for $\alpha \geq 1$). \end{prop}

\begin{proof} By Lemma \ref{exactmean2},
\begin{eqnarray*} \ee_C(c) & = & k - k \sum_{t=0}^{n-1} \frac{B_t (n)_t}{t! k^t}\\ & = & k - k
\sum_{t=0}^{n-1} \frac{B_t n^t}{t! k^t} + k \sum_{t=0}^{n-1}
\frac{B_t (n^t-(n)_t)}{t! k^t}\\ & = & \alpha n - \alpha n
\sum_{t=0}^{n-1} \frac{B_t }{t! \alpha^t} + \alpha n \sum_{t=0}^{n-1}
\frac{B_t (1-\frac{(n)_t}{n^t})}{t! \alpha^t}. \end{eqnarray*} From this and Lemma \ref{estimate0} it follows that
\begin{eqnarray*} & & \left| \ee_C(c) - \left( \alpha n - \alpha n \sum_{t=0}^{\infty} \frac{B_t}{t! \alpha^t}
 + \alpha \sum_{t=0}^{n-1} \frac{B_t {t \choose 2}}{t! \alpha^t}
\right) \right|\\ & \leq & \alpha n \sum_{t=n}^{\infty}
\frac{|B_t|}{t!  \alpha^t} + \frac{\alpha}{2n} \sum_{t=0}^{n-1}
\frac{|B_t| {t \choose 2}^2}{t!  \alpha^t} .\end{eqnarray*} Thus \[
\left| \ee_C(c) - \left( \alpha n - \alpha n \sum_{t=0}^{\infty}
\frac{B_t}{t! \alpha^t} + \alpha \sum_{t=0}^{\infty} \frac{B_t {t
\choose 2}}{t! \alpha^t} \right) \right| \] is at most the ``error
term'' \[ \alpha n \sum_{t=n}^{\infty} \frac{|B_t|}{t!  \alpha^t} +
\alpha \sum_{t=n}^{\infty} \frac{|B_t| {t \choose 2}}{t!  \alpha^t} +
\frac{\alpha}{2n} \sum_{t=1}^{n-1} \frac{|B_t| {t \choose 2}^2}{t!
\alpha^t}. \] From Lemma \ref{estimate1}, the error
term is at most $\frac{C_{\alpha}}{n}$ where $C_{\alpha}$ is a
constant depending on $\alpha$ (and which is independent of $\alpha$
for $\alpha \geq 1$). The result now follows from parts 1 and 2 of
Lemma \ref{bernoulli}. \end{proof}

	Proposition \ref{critvar} estimates the variance of
	$c(\pi)$. The bounds are significantly better than those in
	the literature.

\begin{prop} \label{critvar} Suppose that $k= \alpha n$ with $\alpha > \frac{1}{2
\pi}$. Then\[ \left| Var_{C}(c(\pi)) - n \left( \frac{e^{1/\alpha}
(\alpha^2 e^{2/\alpha} + \alpha^2 - 2 \alpha^2 e^{1/\alpha} -
e^{1/\alpha})}{\alpha^2 (e^{1/\alpha}-1)^4} \right) \right| \leq
A_{\alpha} \] where $A_{\alpha}$ is a constant depending on $\alpha$
(and which is independent of $\alpha$ for $\alpha \geq 1$).
\end{prop}

\begin{proof} From Proposition \ref{exact2} and the expansion of
 partial power sums in \cite{GR}: \[ \sum_{r=0}^{a-1} r^n =
 \frac{a^n}{n+1} \left( a + \sum_{t=0}^{n-1} \frac{B_{t+1}
 (n+1)_{t+1}}{(t+1)! a^t} \right)\] it follows that
\begin{eqnarray*}
& & \ee_C(c^2)\\ & = & k^2 - \frac{n(n+1)}{k^{n-1}} \sum_{j=1}^{k-1}
j^n + \frac{n(nk-n-k)}{k^{n-1}} \sum_{j=1}^{k-1} j^{n-1}\\ & = & k^2 -
nk \left( k+ \sum_{t=0}^{n-1} \frac{B_{t+1} (n+1)_{t+1}}{(t+1)!k^t}
\right)\\ & &  + (nk-n-k) \left( k + \sum_{t=0}^{n-2} \frac{B_{t+1}
(n)_{t+1}}{(t+1)! k^t} \right)\\
& = & -nk - nk  \sum_{t=0}^{n-1} \frac{B_{t+1} (n+1)_{t+1}}{(t+1)!k^t} + (nk-n-k)  \sum_{t=0}^{n-2} \frac{B_{t+1}
(n)_{t+1}}{(t+1)! k^t}\\
& = & -nk - nk \sum_{t=0}^{n-1} \frac{B_{t+1} [(n+1)_{t+1}-(n)_{t+1}]}{(t+1)! k^t} -(n+k) \sum_{t=0}^{n-1} \frac{B_{t+1} (n)_{t+1}}{(t+1)! k^t}\\
& &  - \frac{(nk-n-k) B_n}{k^{n-1}}. \end{eqnarray*} This simplifies to \begin{eqnarray*}
&  & -nk - nk \sum_{t=0}^{n-1} \frac{B_{t+1} (n)_t}{t! k^t} - (nk+k^2) \sum_{t=0}^{n-1} \frac{B_{t+1} (n)_{t+1}}{(t+1)! k^{t+1}}\\ & & - \frac{(nk-n-k) B_n}{k^{n-1}}\\
& = & k^2 - nk \sum_{t=0}^{n-1} \frac{B_{t+1}(n)_t}{t! k^t} - (nk+k^2) \sum_{t=0}^n \frac{B_t (n)_t}{t! k^t} - \frac{(nk-n-k) B_n}{k^{n-1}}\\
& = & \alpha^2 n^2  - \alpha n^2
\sum_{t=0}^{n-1} \frac{B_{t+1}}{t! \alpha^t} 
 - (\alpha n^2 + \alpha^2 n^2)
\sum_{t=0}^n \frac{B_t}{t! \alpha^t}\\
& & + \alpha n^2 \sum_{t=0}^{n-1} \frac{B_{t+1}
(1-\frac{(n)_t}{n^t})}{t! \alpha^t} + (\alpha n^2 + \alpha^2
n^2) \sum_{t=0}^n \frac{B_t (1- \frac{(n)_t}{n^t})}{t! \alpha^t}\\
& & - \frac{(\alpha n^2 - n - \alpha n)B_n}{(\alpha n)^{n-1}}. \end{eqnarray*}

	Lemma \ref{estimate0} implies that the absolute value of the difference between $\ee_C(c^2)$ and  
\begin{eqnarray*}& & 
  \alpha^2 n^2 - \alpha n^2 \sum_{t=0}^{\infty}
   \frac{B_{t+1}}{t! \alpha^t} - (\alpha n^2 + \alpha^2 n^2)
   \sum_{t=0}^{\infty} \frac{B_t}{t! \alpha^t}\\ & & + \alpha n
   \sum_{t=0}^{n-1} \frac{B_{t+1} {t \choose 2}}{t! \alpha^t} +
   (\alpha n + \alpha^2 n) \sum_{t=0}^{n} \frac{B_t {t \choose
   2}}{t! \alpha^t} \end{eqnarray*} is at most
\begin{eqnarray*}
& &  \alpha n^2 \sum_{t=n}^{\infty} \frac{|B_{t+1}|}{t! \alpha^t} +
(\alpha n^2+\alpha^2 n^2) \sum_{t=n+1}^{\infty} \frac{|B_t|}{t!
\alpha^t} + \frac{\alpha}{2} \sum_{t=0}^{n-1} \frac{|B_{t+1}| {t \choose 2}^2}{t! \alpha^t}\\ & & + \frac{(\alpha+\alpha^2)}{2} \sum_{t=0}^{n} \frac{|B_t| {t \choose 2}^2}{t! \alpha^t} + \frac{|\alpha n^2 - n - \alpha n| |B_n|}{(\alpha
n)^{n-1}}. \end{eqnarray*} Thus the difference between $\ee_C(c^2)$ and 
\begin{eqnarray*}& & 
  \alpha^2 n^2 - \alpha n^2 \sum_{t=0}^{\infty}
   \frac{B_{t+1}}{t! \alpha^t} - (\alpha n^2 + \alpha^2 n^2)
   \sum_{t=0}^{\infty} \frac{B_t}{t! \alpha^t}\\ & & + \alpha n
   \sum_{t=0}^{\infty} \frac{B_{t+1} {t \choose 2}}{t! \alpha^t} +
   (\alpha n + \alpha^2 n) \sum_{t=0}^{\infty} \frac{B_t {t \choose
   2}}{t! \alpha^t} \end{eqnarray*} is upper bounded by the ``error term''
\begin{eqnarray*}
& &  \alpha n^2 \sum_{t=n}^{\infty} \frac{|B_{t+1}|\left( 1+ \frac{{t \choose 2}}{n} \right)}{t! \alpha^t} +
(\alpha n^2+\alpha^2 n^2) \sum_{t=n+1}^{\infty} \frac{|B_t| \left( 1+ \frac{{t \choose 2}}{n} \right) }{t!
\alpha^t}\\ & & + \frac{\alpha}{2} \sum_{t=0}^{n-1} \frac{|B_{t+1}| {t \choose 2}^2}{t! \alpha^t} + \frac{(\alpha+\alpha^2)}{2} \sum_{t=0}^{n} \frac{|B_t| {t \choose 2}^2}{t! \alpha^t} + \frac{|\alpha n^2 - n - \alpha n| |B_n|}{(\alpha
n)^{n-1}}. \end{eqnarray*}

	Next it is necessary to bound the five summands in the error
	term. Lemma \ref{estimate1} shows that the
	first four summands in the error term are at most a constant
	depending on $\alpha$ (or a universal constant if $\alpha \geq
	1$).  Since $B_t$ vanishes for $t \geq 3$ odd, $|B_{2t}| \leq
	8 \sqrt{\pi t} (\frac{t}{\pi e})^{2t}$ \cite{Le}, and $2 \pi
	\alpha>1$, the fifth summand in the error term goes to 0 much
	faster than a universal constant.

	The result now follows by combining the above observations
	with Lemma \ref{bernoulli} and Proposition \ref{critmean}. \end{proof}

\subsection{Further properties of the exchangeable pair} \label{further}

	This subsection develops further properties of the
	exchangeable pair $(d,d')$ from Section \ref{exch}. Since we are interested in central limit
	theorems, it is natural to instead study $(W,W')$ where
	$W=\frac{d-\ee_C(d)}{\sqrt{Var_C(d)}}$ and
	$W'=\frac{d'-\ee_C(d)}{\sqrt{Var_C(d)}}$. Note that from Lemma
	\ref{steinmean} one knows $\ee_C(W'-W|\pi)$. In what follows
	we study $\ee_C(W'-W|W)$, which is typically used in normal
	approximation by Stein's method.

\begin{prop} \label{aprox}
\begin{eqnarray*}
& & \ee_C(W'-W|d=r)\\
& = & -\frac{W}{n} + \frac{1}{n \sqrt{Var_C(d)}} \left( \frac{\pp_C(c=r+1)}{\pp_C(d=r)} \frac{(r+1)(n-1)}{n} - \ee_C(d) \right). \end{eqnarray*}
\end{prop}	

\begin{proof} Since $d$ is a function of $\pi$,
\begin{eqnarray*}
\ee_C(W'-W|d=r) & = & \sum_a \frac{a \pp_C(W'-W=a,d=r)}{\pp_C(d=r)}  \\
& = & \sum_a \sum_{\pi: d(\pi)=r}  \frac{a \pp_C(W'-W=a,\pi)}{\pp_C(d=r)} \\
& = & \sum_{\pi: d(\pi)=r} \frac{\pp_C(\pi)}{\pp_C(d=r)}
\sum_a \frac{a \pp_C(W'-W=a,\pi)}{\pp_C(\pi)} \\
& = & \sum_{\pi: d(\pi)=r} \frac{\pp_C(\pi)}{\pp_C(d=r)} \ee_C(W'-W|\pi). \end{eqnarray*}
	
	By Lemma \ref{steinmean}, this is equal to
\begin{eqnarray*}
& & \frac{1}{\sqrt{Var_C(d)}} \left( \sum_{\pi: d(\pi)=r} \frac{\pp_C(\pi)}{\pp_C(d=r)} \left[ -\frac{r}{n} + \frac{n-1}{n} \ii_{\chi_n(\pi)=1} \right] \right)\\
& = &  \frac{ -r + \ee_C(d)}{n \sqrt{Var_C(d)}} +  \frac{ \left( \sum_{\pi: d(\pi)=r} \frac{\pp_C(\pi)}{\pp_C(d=r)} \left[  (n-1) \ii_{\chi_n(\pi)=1} - \ee_C(d) \right] \right)  }{n \sqrt{Var_C(d)}}\\
& = & -\frac{W}{n} + \frac{1}{n \sqrt{Var_C(d)}} \left( \frac{(n-1) \sum_{\pi: d(\pi)=r} \pp_C(\pi) \ii_{\chi_n(\pi)=1}}{\pp_C(d=r)} - \ee_C(d) \right)\\
& = & -\frac{W}{n} + \frac{1}{n \sqrt{Var_C(d)}} \left( \frac{\pp_C(c=r+1)}{\pp_C(d=r)} \frac{(r+1)(n-1)}{n} - \ee_C(d) \right).
\end{eqnarray*}
\end{proof}

	In most examples of Stein's method for normal approximation of
	a random variable $W$, there is an exchangeable pair $(W,W')$
	such that $\ee(W'|W) = (1-\lambda)W$. There are two recent
	papers (\cite{RR},\cite{C}) in which the Stein technique has
	been extended to handle the case where $\ee(W'|W) = (1 -
	\lambda)W + G(W)$ where $G(W)$ is small. The bounds in these
	papers require that $\frac{\ee(|G(W)|)}{\lambda}$ goes to
	0. Proposition \ref{nogood} uses interesting properties (Lemma \ref{newton}) and asymptotics of
	Eulerian numbers to prove that for our exchangeable pair
	$\frac{\ee(|G(W)|)}{\lambda}$ is bounded away from 0, even for
	$C_{\infty,k}$ (the uniform distribution on the symmetric
	group), where we know that a central limit theorem holds.

\begin{lemma} \label{newton} Let $A_{n,k}$ denote the number of permutations on $n$ symbols with $k-1$ descents. 
\begin{enumerate}
\item If $n$ is odd and $0 \leq r \leq n-1$, then $(r+1) A_{n-1,r+1} \geq (n-r) A_{n-1,r}$ if and only if $0 \leq r \leq \frac{n-1}{2}$.
\item If $n$ is even and $0 \leq r \leq n-1$, then $(r+1) A_{n-1,r+1} \geq (n-r) A_{n-1,r}$ if and only if $0 \leq r \leq \frac{n}{2}-1$.
\end{enumerate}
\end{lemma}

\begin{proof} Suppose first that $n$ is odd. For the if part, we proceed by
 reverse induction on $r$; thus the base case is $r=\frac{n-1}{2}$,
 and then the inequality is an equality since $A_{n-1,k}=A_{n-1,n-k}$ for all $k$. A result of Frobenius
 \cite{Fr} is that the polynomial $\sum_{k \geq 1} z^k A_{n-1,k}$ has
 only real roots. Thus an inequality of Newton (page 52 of \cite{HLP})
 implies that \[ \frac{(r+1) A_{n-1,r+1}}{(n-r) A_{n-1,r}} \geq
 \frac{A_{n-1,r+2}}{A_{n-1,r+1}} \frac{(r+2)(n-r-1)}{(n-r-2)(n-r)}.\] By
 the induction hypothesis the right hand side is at least
 $\frac{(n-r-1)^2}{(n-r-2)(n-r)} > 1$. The only if part follows from the if part since $A_{n-1,k}=A_{n-1,n-k}$ for all $k$.

	The case of $n$ even is similar. For the if part, we proceed
	by reverse induction on $r$. The induction step is the same
	but base case $r=\frac{n}{2}-1$ is not automatic. However it
	follows using Newton's inequality of the previous paragraph
	together with the symmetry property $A_{n-1,k}=A_{n-1,n-k}$:
	\begin{eqnarray*} (A_{n-1,\frac{n}{2}})^2 & \geq &
	(A_{n-1,\frac{n}{2}-1}) (A_{n-1,\frac{n}{2}+1})
	\frac{(\frac{n}{2}+1)}{(\frac{n}{2}-1)}\\ & = &
	(A_{n-1,\frac{n}{2}-1})^2 \frac{(\frac{n}{2}+1)}{(\frac{n}{2}-1)}\\ & >
	& (A_{n-1,\frac{n}{2}-1})^2
	(1+\frac{2}{n})^2. \end{eqnarray*} Now take square roots.
	The only if part follows from the if part since
	$A_{n-1,k}=A_{n-1,n-k}$ for all $k$. \end{proof}

\begin{prop} \label{nogood} Let
 $U_n$ denote the uniform distribution $U_n$ on $S_n$. Let $(W,W')$ be the exchangeable pair of this subsection,
 so that by Proposition \ref{aprox}, $\ee(W'|W) = (1-\lambda)W + G(W)$
 with $\lambda=\frac{1}{n}$. Then $\frac{\ee_{U_n}(|G(W)|)}{\lambda}$ is
 bounded away from 0 as $n \rightarrow \infty$. \end{prop}

\begin{proof} It is elementary that $\ee_{U_n}(d)=\frac{n-1}{2}$. Thus Proposition \ref{aprox} implies that 
\[ \ee_{U_n}(|G(W)|) =  \frac{(n-1)}{n^2 \sqrt{Var_{U_n}(d)}} \sum_{r=0}^{n-1} \left| (r+1) \pp_{U_n}(c=r+1) - \frac{n}{2} \pp_{U_n}(d=r) \right| .\] Using the equality \[ \pp_{U_n}(d=r) = \frac{(r+1)}{n} \pp_{U_n}(c=r+1) + \frac{(n-r)}{n} \pp_{U_n}(c=r), \] this simplifies to \[ \frac{(n-1)}{2 n^2 \sqrt{Var_{U_n}(d)}} \sum_{r=0}^{n-1} \left| (r+1) \pp_{U_n}(c=r+1) - (n-r)
\pp_{U_n}(c=r) \right| .\] By Lemma \ref{ful}, this further
simplifies to \[ \frac{(n-1)}{ (n-1)! 2 n^2 \sqrt{Var_{U_n}(d)}}
\sum_{r=0}^{n-1} \left| (r+1) A_{n-1,r+1} - (n-r) A_{n-1,r} \right|,
\] where $A_{n,k}$ denotes the number of permutations in $S_n$ with
$k-1$ descents.

	Now suppose that $n$ is odd. Then the previous paragraph,
Lemma \ref{newton}, and the symmetry $A_{n-1,k}=A_{n-1,n-k}$ imply
that  $\ee_{U_n}(|G(W)|)$ is equal to \begin{eqnarray*} &  & \frac{(n-1)}{ (n-1)!
n^2 \sqrt{Var_{U_n}(d)}} \sum_{r=0}^{\frac{n-3}{2}} \left( (r+1)A_{n-1,r+1} -
(n-r) A_{n-1,r} \right) \\ & = & \frac{(n-1)}{(n-1)! n^2 \sqrt{Var_{U_n}(d)}}
\left( \sum_{r=1}^{\frac{n-1}{2}} r A_{n-1,r} -
\sum_{r=1}^{\frac{n-3}{2}} (n-r) A_{n-1,r} \right)\\ & = &
\frac{(n-1)}{(n-1)! n^2 \sqrt{Var_{U_n}(d)}} \left( \frac{n+1}{2}
A_{n-1,\frac{n-1}{2}} + \sum_{r=1}^{\frac{n-1}{2}} (2r-n) A_{n-1,r}
\right) \\ & = & \frac{(n-1)}{n^2 \sqrt{Var_{U_n}(d)}} \left(
\frac{n+1}{2} \frac{A_{n-1,\frac{n-1}{2}}}{(n-1)!} -\ee_{U_{n-1}}
|d-\ee_{U_{n-1}}(d)| \right) \\ & \geq & \frac{(n-1)}{n^2
\sqrt{Var_{U_n}(d)}} \left( \frac{n+1}{2} \frac{A_{n-1,\frac{n-1}{2}}}{2
(n-1)!} - \sqrt{Var_{U_{n-1}}(d)} \right). \end{eqnarray*} A similar
argument for $n$ even shows that \[ \ee_{U_n}(|G(W)|) \geq
\frac{(n-1)}{n^2 \sqrt{Var_{U_n}(d)}} \left( \frac{ n
A_{n-1,\frac{n}{2}}}{2 (n-1)!} - \sqrt{Var_{U_{n-1}}(d)} \right). \]

	To conclude the argument, note that for $n \geq 2$,
	$Var_{U_n}(d) = \frac{n+1}{12}$, and that $\frac{A_{n-1,\lceil
	\frac{n-1}{2} \rceil}}{(n-1)!}$ is asymptotic to
	$\sqrt{\frac{6}{n \pi}}$ 
	\cite{CKSS}, \cite{T}. Thus $n \ee_{U_n}(|G(W)|)$ is bounded away from
	0, as desired. \end{proof}

\subsection{A Sampling Algorithm} \label{sample}

	To conclude the paper, we record a fast algorithm for drawing
	samples from the measure $R_{k,n}$. Since $C_{k,n}$ can be
	obtained by sampling from $R_{k,n}$ and then performing
	uniformly at random one of the n cyclic rotations of the
	bottom row in the two line form of a permutation, we only give
	an algorithm for sampling from $R_{k,n}$. This algorithm
	should be quite useful for empirically studying the transition
	from Poisson to normal behavior.

	We use the terminology that if $\pi$ is a permutation on $n$ symbols, the permutation $\tau$ on $n+1$ symbols obtained by inserting $n+1$ after position $j$ (with $0 \leq j \leq n$) is defined by $\tau(i)=\pi(i)$ for $1 \leq i \leq j$, $\tau(j+1) = n+1$, and $\tau(i)=\pi(i-1)$ for $j+2 \leq i \leq n+1$. For instance inserting 5 after position 2 in the permutation $3 \ 4 \ 1 \ 2 $ gives the permutation $3 \ 4 \ 5 \ 1 \ 2$.

\begin{prop} Starting with the identity permutation in $S_1$, transition from an element $\pi$ of $S_{n}$ to an element $\tau$ of $S_{n+1}$ by inserting $n+1$ as described by the following 2 cases:
\begin{enumerate}
\item If either $j=n$ or $\pi(j)> \pi(j+1)$ and $1 \leq j \leq n-1$, the chance of inserting $n+1$ after $j$ is $\frac{n+k-d(\pi)}{k(n+1)}$.
\item If either $j=0$ or $\pi(j)< \pi(j+1)$ and $1 \leq j \leq n-1$, the chance of inserting $n+1$ after $j$ is $\frac{k-d(\pi)-1}{k(n+1)}$.
\end{enumerate} Then after running the algorithm for $n-1$ steps, the distribution obtained on $S_n$ is precisely $R_{k,n}$.
\end{prop}

\begin{proof} First note that the transition probabilities of the algorithm sum to 1, since \[ (d(\pi)+1) \left( \frac{n+k-d(\pi)}{k(n+1)} \right) + (n-d(\pi)) \left( \frac{k-d(\pi)-1}{k(n+1)} \right) =1.\] Now observe that if $\tau$ is obtained from $\pi$ by a Case 1 move, then $d(\tau)=d(\pi)$. Thus from the formula for $R_{k,n}$, $\frac{R_{k,n+1}(\tau)}{R_{k,n}(\pi)} = \frac{n+k-d(\pi)}{k(n+1)}$. Similarly if $\tau$ is obtained from $\pi$ by a Case 2 move, then $d(\tau)=d(\pi)+1$. Thus $\frac{R_{k,n+1}(\tau)}{R_{k,n}(\pi)} = \frac{k-d(\pi)-1}{k(n+1)}$. \end{proof}

\section{Acknowledgements} The author was partially supported by NSA grant MDA904-03-1-0049.

\end{document}